\documentclass{article}
\usepackage{amssymb,latexsym}
\usepackage{graphicx}
\usepackage{amsmath}
\usepackage{amsthm}
\usepackage{amscd}
\usepackage{bm}
\begin{document}
\newcommand{\D}{\displaystyle}
\newcommand{\DF}[2]{\frac{\D#1}{\D#2}}
\def\Xint#1{\mathchoice
{\XXint\displaystyle\textstyle{#1}}%
{\XXint\textstyle\scriptstyle{#1}}%
{\XXint\scriptstyle\scriptscriptstyle{#1}}%
{\XXint\scriptscriptstyle\scriptscriptstyle{#1}}%
\!\int}
\def\XXint#1#2#3{{\setbox0=\hbox{$#1{#2#3}{\int}$}
\vcenter{\hbox{$#2#3$}}\kern-.5\wd0}}
\def\ddashint{\Xint=}
\def\dashint{\Xint-}
\pagestyle{plain}
\newpage
\normalsize
\title{A Theorem on Frequency Function for Multiple-Valued Dirichlet Minimizing Functions
\thanks{Department of mathematics, Rice University, Houston, TX 77005, U.S.A; weizhu@math.rice.edu}}
\author{Wei Zhu}
\maketitle
\begin{abstract}
This paper discusses the frequency function of multiple-valued
Dirichlet minimizing functions in the special case when the domain
and range are both two dimensional. It shows that the frequency
function must be of value $k/2$ for some nonnegative integer $k$.
Futhermore, by looking at the blowing-up functions, we
characterize the local behavior of the original Dirichlet
minimizing function.
\end{abstract}
\newpage
\newtheorem{theorem}{Theorem}[section]
\newtheorem{definition}{Definition}[section]
\newtheorem{remark}{Remark}[section]
\newtheorem{corollary}{Corollary}[section]
\newtheorem{proposition}{Proposition}[section]
\newtheorem{lemma}{Lemma}[section]
\section{Introduction}
Frequency function for multiple-valued functions was introduced by
Almgren in $\cite{af}$ to study the branching behavior for
multiple-valued Dirichlet minimizing functions:
$$N(r)=\DF{r\int_{\mathbb{B}_r^m(a)}|Df|^2}{\int_{\partial
\mathbb{B}_r^m(a)}|f|^2}.$$ For Dirichlet minimizing functions,
$N(r)$ is nondecreasing in $r$. Almgren establishes this
monotonicity by certain range and domain deformations, called
``squashing'' and ``squeezing''. The monotonicity property enables
one to prove by dimension reduction that such multiple-valued
functions have branched sets of codimension at least two. \\
To get a better idea about this frequency function, consider a
harmonic function on $\mathbb{R}^2$, and express it in terms of
polar coordinates: $u(r,\theta)$. If we fix $r$, we can expand the
resulting function of $\theta$ as a Fourier series. Now as $r$
decreases, the higher frequency terms in the Fourier series die
off faster than the lower frequency terms.\\
Monotonicity of frequency functions have been used in some other
work, see $\cite{gl},\cite{lf},\cite{gs}$.\\
The motivation of this paper was trying to characterize
multiple-valued Dirichlet minimizing functions
$f:\mathbb{R}^2\to\mathbb{Q}_2(\mathbb{R}^2)$, which is
homogeneous of some positive degree. There are a lot of them like
$z^{1/2},\pm z,z^{3/2}$. In general, any function of the form
$$z^N,\;\mbox{for some positive real number}\;N$$
could be a candidate. One thing worth mentioning is that the
frequency function at the origin of $z^N$ is exactly $N$. However,
not every $N$ gives a 2-valued function because the function has
to match up itself once going around the circle one time. For
example, consider the function
$$f:(r,\theta)\to
[[(r^N\cos(N\theta),r^N\sin(N\theta))]]+[[(-r^N\cos(N\theta),-r^N\sin(N\theta))]],$$
when $N=1/4$.
$$f(r,0)=[[(r^N,0)]]+[[(-r^N,0)]],f(r,2\pi)=[[(0,r^N)]]+[[(0,-r^N)]].$$
They do not match. We will see that in this case, only by choosing
$N=k/2$ for some positive integer $k$ makes $f$ a well-defined
2-valued function. This is basically the main ingredient of the
proof of our main theorem, matching up values for $\theta=0$ and
$\theta=2\pi$. More precisely, given a multiple-valued Dirichlet
minimizing function $f:\mathbb{R}^2\to\mathbb{Q}_2(\mathbb{R}^2)$,
with $\mathcal{N}(0)=N$, we use the blowing-up analysis to get a
Dirichlet minimizing function
$g:\mathbb{R}^2\to\mathbb{Q}_2(\mathbb{R}^2)$ of homogeneous
degree $N$ with the same frequency $N$ at the origin. By doing the
matching up business for $g$, we succeed in proving $N=k/2$ for
some nonnegative integer $k$. A by-product of this proof is the
characterization of local behavior of original function $f$ near
the origin.\\
So a natural question is whether we have similar results in higher
dimension, either domain or the range, or higher multiplicity $Q$.
For functions $f:\mathbb{R}^2\to\mathbb{Q}_3(\mathbb{R}^2)$, it
would take a lot more work to matching up values. Therefore some
other easier methods are expected to give a full answer to this
question in general.
\section{Preliminaries}
We refer to $\cite{af},\cite{zw1}$ for most of notations,
definitions and known results about multiple-valued functions. For
reader's convenience, here we state some useful results. The
proofs of them can be found in $\cite{af}$.
\begin{theorem}[$\cite{af}, \S 2.6$]
{\bfseries Hypotheses:}\\
(a) $0<r_0<\infty$.\\
(b) $A\subset \mathbb{R}^m$ is connected, open, and bounded with
$\mathbb{U}_{r_0}^m(0)\subset A$. $\partial A$
is an $m-1$ dimensional submanifold of $\mathbb{R}^m$ of class 1.\\
(c) $f:A\rightarrow \mathbb{Q}$ is strictly defined and is Dir minimizing.\\
(d) $D,H,N:(0,r_0)\rightarrow\mathbb{R}$ are defined for $0<r<r_0$
by setting
\begin{equation*}
\begin{split}
D(r)&=\mbox{Dir}(f;\mathbb{B}_r^m(0))\\
H(r)&=\int_{\partial \mathbb{B}_r^m(0)}\mathcal{G}(f(x),Q[[0]])^2
d\mathcal{H}^{m-1} x\\
N(r)&=rD(r)/H(r)\;\mbox{provided}\;H(r)>0.
\end{split}
\end{equation*}
(e) $\mathcal{N}:A\rightarrow \mathbb{R}$ is defined for $x\in A$ by
setting
$$\mathcal{N}(x)=\lim_{r\downarrow 0} r
Dir(f;\mathbb{B}_r^m(x))/\int_{\partial
\mathbb{B}_r^m(x)}\mathcal{G}(f(z),Q[[0]])^2d\mathcal{H}^{m-1}z$$
provided
this limit exists.\\
(f) $H(r)>0$ for some $0<r<r_0$.\\
{\bfseries Conclusions.}\\
(1) $\eta\circ f\in\mathcal{Y}_2(A,\mathbb{R}^n)$ is Dir minimizing
and harmonic.\\
(2) $N(r)$ is defined for each $0<r<r_0$ and is nondecreasing.\\
(3) $\mathcal{N}(0)=\lim_{r\downarrow 0}N(r)$ exists.\\
(4) \\$A=A\cap \{x:\mbox{for some}\;0<r<\mbox{dist}(x,\partial
A),\int_{\partial \mathbb{U}_r^m(x)}\mathcal{G}(f(z),Q[[0]])^2
d\mathcal{H}^{m-1}
>0\}.$\\
(5) $\mathcal{N}(x)$ is well defined for each $x\in A$ and is upper
semi-continuous as a function of $x$.\\
(6) In case $N(r)=\mathcal{N}(0)$ for $\mathcal{L}^1$ almost all
$0<r<r_0$, then
$$f(x)=\mathbb{\mu}[(r/r_1)^{\mathcal{N}(0)}]_\sharp f(r_1 x/|x|)$$
for $\mathcal{L}^{m-1}$ almost all $x\in \partial
\mathbb{U}_r^m(0)$ and each $0<r_1<r_0$.
\end{theorem}
\begin{theorem}[$\cite{af}, \S 2.13$]
{\bfseries Hypotheses.}\\
(a) In case $m=2$, $\omega_{2.13}=1/Q$.\\
(b) In case $m\ge 3, 0<\epsilon_Q<1$ is as defined as in
$\cite{af},\S 2.11$ and $0<\omega_{2.13}<1$ is defined by the
requirement
$$m-2+2\omega_{2.13}=(m-2)(1+\epsilon_Q)/(1-\epsilon_Q).$$
(c)
$$\Gamma_{2.13}=4^{1-\omega_{2.13}}[2^{m/2}/(1-2^{-\omega_{2.13}})+3\cdot
2^{m-1+\omega_{2.13}}](m\alpha(m))^{-1/2}
\mbox{Lip}(\xi)^2\mbox{Lip}(\xi^{-1}).$$ (d)
$f\in\mathcal{Y}_2(\mathbb{R}^m,\mathbb{Q})$ is strictly defined
and $f|\mathbb{U}_1^m(0)$ is Dir minimizing with
Dir$(f;\mathbb{B}_1^m(0))>0$.\\
(e) $\mathcal{N}(0)=\lim_{r\downarrow 0} r
Dir(f;\mathbb{B}_r^m(0))/\int_{x\in\partial
\mathbb{B}_r^m(0)}\mathcal{G}(f(x),Q[[0]])^2 d\mathcal{H}^{m-1} x$.\\
{\bfseries Conclusions.}\\
(1) For each $z\in \mathbb{U}_1^m(0),0<r<1-|z|$, and $0<s\le 1$,
$$\mbox{Dir}(f;\mathbb{B}_{sr}^m(z))\le s^{m-2+2\omega_{2.13}}\mbox{Dir}(f;\mathbb{B}_r^m(z)).$$
(2) Whenever $0<\delta<1$ and $p,q\in \mathbb{B}_{1-\delta}^m(0)$,
$$\mathcal{G}(f(p),f(q))\le \Gamma_{2.13}\delta^{-m/2}\mbox{Dir}(f;\mathbb{B}_1^m(0))^{1/2}|p-q|^{\omega_{2.13}},$$
in particular, $f|\mathbb{B}_{1-\delta}^m(0)$ is H$\ddot{o}$lder continuous with exponent $\omega_{2.13}$. \\
(3) Either $f(0)=Q[[0]]$ and $\mathcal{N}(0)\ge \omega_{2.13}$ or
$f(0)\not= Q[[0]]$ and $\mathcal{N}(0)=0$.\\
(4) Suppose $f(0)=Q[[0]]$ and
$1/2>r(1)>r(2)>r(3)>\cdot\cdot\cdot>0$ with $0=\lim_{i\rightarrow
\infty} r(i)$. Then there is a subsequence
$i_1,i_2,i_3,\cdot\cdot\cdot$ of $1,2,3,\cdot\cdot\cdot$ and a
function $g:\mathbb{B}_1^m(0)\rightarrow \mathbb{Q}$ with the
following
properties:\\
(a) $g$ is the uniform limit as $k\rightarrow \infty$ of the
functions
$$\mu(\mbox{Dir}(f\circ\mu[r(i_k)];\mathbb{B}_1^m(0))^{-1/2})_{\sharp}\circ f\circ\mu[r(i_k)]|\mathbb{B}_1^m(0).$$
(b) $g|\mathbb{U}_1^m(0)\in\mathcal{Y}_2(\mathbb{U}_1^m(0),\mathbb{Q})$ is Dir minimizing with Dir($g;\mathbb{B}_1^m(0))=1$.\\
(c) $\int_{x\in\partial \mathbb{B}_1^m(0)}\mathcal{G}(g(x),Q[[0]])^2 d\mathcal{H}^{m-1} x=1/\mathcal{N}(0)$.\\
(d) $g(0)=Q[[0]]$ and for each $x\in \mathbb{B}_1^m(0)\sim\{0\}$,
$$g(x)=\mu(|x|^{\mathcal{N}(0)})_{\sharp}\circ g(x/|x|).$$
(e) For each $p,q\in \mathbb{B}_1^m(0)$,
$$\mathcal{G}(g(p),g(q))\le
2^{m/2-\omega_{2.13}+\mathcal{N}(0)}\Gamma_{2.13}|p-q|^{\omega_{2.13}}.$$
(5) Corresponding to each bounded open set $A$ such that $\partial
A$ is a compact $m-1$ dimensional submanifold of $\mathbb{R}^m$ of
class 1, there is a constant $0<\Gamma_A<\infty$ with the following
property. Whenever $g\in\mathcal{Y}_2(A,\mathbb{Q})$ is Dir
minimizing and $p,q\in A$,
$$\mathcal{G}(g(p),g(q))\le \Gamma_A\mbox{Dir}(g;A)^{1/2}\sup\{\mbox{dist}(p,\partial A)^{-m/2},\mbox{dist}(q,\partial A)^{-m/2}\}
|p-q|^{\omega_{2.13}}.$$
\end{theorem}
\begin{theorem}[$\cite{af},\S 2.14$]
(1) Let $\mu\in\{1,2,\cdot\cdot\cdot,Q\}$ and suppose
$f_1,f_2,\cdot\cdot\cdot,f_Q\in\mathcal{Y}_2(\mathbb{U}_1^\mu(0),\mathbb{R}^n)$
are strictly defined. Then $f=\sum_{i=1}^Q [[f_i]]\in
\mathcal{Y}_2(\mathbb{U}_1^\mu,\mathbb{Q})$. Furthermore, in case
$f$ is Dir
minimizing, so is each $f_i,i=1,2,\cdot\cdot\cdot,Q$.\\
(2)Suppose $f\in\mathcal{Y}_2(\mathbb{U}_1^m(0),\mathbb{Q})$ is
strictly defined and Dir minimizing. Then the function
$$\sigma:\mathbb{U}_1^m(0)\rightarrow \{1,2,\cdot\cdot\cdot,Q\},$$
$$\sigma(x)=\mbox{card}\;[\mbox{spt}(f(x))]\;\mbox{for}\;x\in
\mathbb{U}_1^m(0),$$ is lower semi-continuous, the set
$$\Sigma=\mathbb{U}_1^m(0)\cap\{x:\sigma\;\mbox{is not continuous at}\;x\}$$
is closed in $\mathbb{U}_1^m(0)$ with Hausdorff dimension not
exceeding $m-2$, and the set $\mathbb{U}_1^m(0)\sim \Sigma$ is
open and path connected. Furthermore, there exist
$J\in\{1,2,\cdot\cdot\cdot,Q\}$ and
$k_1,k_2,\cdot\cdot\cdot,k_J\in\{1,2,\cdot\cdot\cdot,Q\}$ with
$k_1+k_2+\cdot\cdot\cdot+k_J=Q$ with the following properties:
whenever $W\subset \mathbb{U}_1^m(0)\sim \Sigma$ is open and
simply connected, there are harmonic functions
$f_1,f_2,\cdot\cdot\cdot,f_J:W\rightarrow \mathbb{R}^n$ such that
$f(x)=\sum_{i=1}^J k_i[[f_i(x)]]$ and
$J=\mbox{card}\;\{f_1(x),f_2(x),\cdot\cdot\cdot,f_J(x)\}$ for each
$x\in W$.
\end{theorem}
\section{Main Theorem}
\begin{theorem}
{\bfseries Hypotheses.}\\
(a) $m=2,n=2,Q=2,\omega_{2.13}=1/2$.\\
(b) $f\in\mathcal{Y}_2(\mathbb{R}^2,\mathbb{Q}(\mathbb{R}^2))$ is
strictly defined and $f|\mathbb{U}_1^2(0)$ is Dir minimizing with
Dir($f;\mathbb{B}_1^2(0))>0$.\\
(c) $f(0)=2[[0]]$,\\
(d) $\mathcal{N}(0)=\lim_{r\downarrow 0}
r\mbox{Dir}(f;\mathbb{B}_r^2(0))/\int_{x\in\partial
\mathbb{B}_r^2(0)}
\mathcal{G}(f(x),2[[0]])^2d\mathcal{H}^1x$\\
{\bfseries Conclusion.}
$$\mathcal{N}(0)=k/2,\;\mbox{for some positive integer}\;k.$$
\end{theorem}
\begin{proof}
According to Theorem 2.2, we know that $\mathcal{N}(0)\ge
\omega_{2.13}>0$ and suppose $1/2> r(1)> r(2)>
r(3)>\cdot\cdot\cdot> 0$ with $0=\lim_{i\rightarrow \infty} r(i)$,
then there is a subsequence
$i_1,i_2,i_3,\cdot\cdot\cdot$ of $1,2,3,\cdot\cdot\cdot$ and a function $g:\mathbb{B}_1^2(0)\rightarrow \mathbb{Q}$ with the following properties:\\
(1) $g$ is the uniform limit as $k\rightarrow \infty$ of the
functions
$$\mu(\mbox{Dir}(f\circ\mu[r(i_k)];\mathbb{B}_1^2(0))^{-1/2})_{\sharp}\circ f\circ\mu[r(i_k)]|\mathbb{B}_1^2(0).$$
(2) $g|\mathbb{U}_1^2(0)\in\mathcal{Y}_2(\mathbb{U}_1^2(0),\mathbb{Q})$ is Dir minimizing with Dir($g;\mathbb{B}_1^2(0))=1$.\\
(3) $\int_{x\in\partial \mathbb{B}_1^2(0)}\mathcal{G}(g(x),2[[0]])^2 d\mathcal{H}^1 x=1/\mathcal{N}(0)$.\\
(4) $g(0)=2[[0]]$ and for each $x\in \mathbb{B}_1^2(0)\sim\{0\}$,
$$g(x)=\mu(|x|^{\mathcal{N}(0)})_{\sharp}\circ g(x/|x|).$$
(5) For each $p,q\in \mathbb{B}_1^2(0)$,
$$\mathcal{G}(g(p),g(q))\le
2^{1/2+\mathcal{N}(0)}\Gamma_{2.13}|p-q|^{1/2}.$$ First of all, we
claim
$$\Sigma(g)=\emptyset,\;\mbox{or}\;\{0\}.$$ This comes from the fact
that $g$ is homogeneous of degree $\mathcal{N}(0)$. If $\sigma$ is
not continuous at some nonzero point $y$, $\sigma$ is not continuous
at every point on the ray $ty,t\in (0,1)$. Then the Hausdorff
dimension of $\Sigma$ is at least one, which is in contradiction to
Theorem 2.3(2). \\
Rest of the proof is divided into two cases: $\Sigma=\emptyset$
and $\Sigma=\{0\}.$
\subsection{$\Sigma=\{0\}$}
If $\Sigma=\{0\}$, applying Theorem 2.3(2) to the function $g$, we
get $J=2$. This is because otherwise if $J=1$, then for any point
$x\in \mathbb{U}_1^2(0)\sim \{0\}$, $\sigma(x)=1$. Therefore
$\sigma$ is a constant function on $\mathbb{U}_1^2(0)$, which
means $\Sigma=\emptyset$, a
contradiction to our assumption. \\
Take $W=\mathbb{U}_1^2(0)\sim \{(x,0),x\ge 0\}$ in Theorem 2.3(2),
we have
$$g(x)=\sum_{i=1}^2 [[h_i(x)]],x\in W$$
for harmonic functions $h_i:W\rightarrow \mathbb{R}^2,i=1,2$ and $h_1(x)\not= h_2(x),\forall x\in W$.\\
For simplicity, we denote $\mathcal{N}(0)$ as $N$. Since $g$ is
homogeneous of degree $N$, so is $h_i,i=1,2$. Hence we can write
$$g(r,\theta)=[[r^Ng_1(\theta)]]+[[r^Ng_2(\theta)]],0< r\le 1,0<
\theta< 2\pi$$
where $g_i:(0,2\pi)\rightarrow \mathbb{R}^2,i=1,2$, and $r^N g_i$ is harmonic, $i=1,2$.\\
Moreover, in spirit of Theorem 2.3(1), $r^N g_i$ must be Dir
minimizing, hence conformal on $W$ for $i=1,2$.\\
Let
$$g_1(\theta)=(g_1^1(\theta),g_1^2(\theta)),g_2(\theta)=(g_2^1(\theta),g_2^2(\theta)).$$
The Laplacian operator in polar coordinate can be expressed as
$$\Delta=\DF{\partial^2}{\partial r^2}+\DF{1}{r}\DF{\partial}{\partial r}+\DF{1}{r^2}\DF{\partial^2}{\partial \theta^2}.$$
Do the computation, we have
$$\DF{\partial^2}{\partial r^2}(r^Ng_i^j(\theta))=N(N-1)r^{N-2}g_i^j(\theta),$$
$$\DF{1}{r}\DF{\partial}{\partial r}(r^Ng_i^j(\theta))=Nr^{N-2}g_i^j(\theta),$$
$$\DF{1}{r^2}\DF{\partial^2}{\partial\theta^2}(r^Ng_i^j(\theta))=r^{N-2}[g_i^j(\theta)]''.$$
Therefore, $\Delta(r^Ng_i^j(\theta))=r^{N-2}[N^2g_i^j(\theta)+[g_i^j(\theta)]'']=0,i=1,2,j=1,2$.\\
Hence we can represent $g_i^j(\theta)$ as
$$g_i^j(\theta)=a_i^j\cos(N\theta)+b_i^j\sin(N\theta),i=1,2,j=1,2,\;\mbox{for some contants}\;a_i^j,b_i^j.$$
Denote
$$g(r,\theta)=r^N[[(a\cos(N\theta)+b\sin(N\theta),c\cos(N\theta)+d\sin(N\theta))]]$$
$$+r^N[[(\tilde{a}\cos(N\theta)+
\tilde{b}\sin(N\theta),\tilde{c}\cos(N\theta)+\tilde{d}\sin(N\theta))]],0<r<1,0<\theta<2\pi$$
$\mbox{where}\;a,b,c,d,\tilde{a},
\tilde{b},\tilde{c},\tilde{d}\;\mbox{are constants}.$\\
Denote \\$h_1=r^Ng_1(\theta)=(r^N(a\cos(N\theta)+b\sin(N\theta)),r^N(c\cos(N\theta)+d\sin(N\theta)))=(f_1,f_2)$. \\
In the polar coordinate,
$$\DF{\partial}{\partial x}=\DF{\partial}{\partial r}\cdot \cos\theta+\DF{\partial}{\partial \theta}\cdot
(\DF{-\sin\theta}{r}),$$
$$\DF{\partial}{\partial y}=\DF{\partial}{\partial r}\cdot \sin\theta+\DF{\partial}{\partial \theta}\cdot (\DF{\cos\theta}{r}).$$
Do the computation,
$$\DF{\partial f_1}{\partial x}=\DF{\partial f_1}{\partial r}\cos\theta+\DF{\partial f_1}{\partial \theta}(\DF{-\sin\theta}{r})$$
\begin{equation*}
\begin{split}
&=Nr^{N-1}(a\cos(N\theta)+b\sin(N\theta))\cos\theta+r^N(-aN\sin(N\theta)+bN\cos(N\theta))\DF{-\sin\theta}{r}\\
&=Nr^{N-1}(a\cos(N\theta)\cos\theta+b\sin(N\theta)\cos\theta+a\sin(N\theta)\sin\theta-b\cos(N\theta)\sin\theta)\\
&=Nr^{N-1}(a\cos((N-1)\theta)+b\sin((N-1)\theta)).
\end{split}
\end{equation*}
Similarly, $\DF{\partial f_2}{\partial
x}=Nr^{N-1}(c\cos((N-1)\theta)+d\sin((N-1)\theta))$.\\
$$\DF{\partial f_1}{\partial y}=\DF{\partial f_1}{\partial r}\sin\theta+\DF{\partial f_1}{\partial \theta}(\DF{\cos\theta}{r})$$
\begin{equation*}
\begin{split}
&=Nr^{N-1}(a\cos(N\theta)+b\sin(N\theta))\sin\theta+r^N(-aN\sin(N\theta)+bN\cos(N\theta))\DF{\cos\theta}{r}\\
&=Nr^{N-1}(a\cos(N\theta)\sin\theta+b\sin(N\theta)\sin\theta-a\sin(N\theta)\cos\theta+b\cos(N\theta)\cos\theta)\\
&=Nr^{N-1}(b\cos((N-1)\theta)-a\sin((N-1)\theta))
\end{split}
\end{equation*}
Similarly, $\DF{\partial f_2}{\partial y}=Nr^{N-1}(d\cos((N-1)\theta)-c\sin((N-1)\theta))$.\\
Let $(N-1)\theta=\phi$,\\
$$\DF{\partial h_1}{\partial x}=(\DF{\partial f_1}{\partial x},\DF{\partial f_2}{\partial x})=Nr^{N-1}(a\cos\phi+b\sin\phi,c\cos\phi+d\sin\phi),$$
$$\DF{\partial h_1}{\partial y}=(\DF{\partial f_1}{\partial y},\DF{\partial f_2}{\partial y})
=Nr^{N-1}(b\cos\phi-a\sin\phi,d\cos\phi-c\sin\phi).$$ $|\DF{\partial
h_1}{\partial
x}|^2=N^2r^{2(N-1)}[a^2\cos^2\phi+b^2\sin^2\phi+2ab\sin\phi\cos\phi+c^2\cos^2\phi+d^2\sin^2\phi+2cd\sin\phi
\cos\phi]$\\
$|\DF{\partial h_1}{\partial
y}|^2=N^2r^{2(N-1)}[b^2\cos^2\phi+a\sin^2\phi-2ab\sin\phi\cos\phi+
d^2\cos^2\phi+c^2\sin^2\phi-2cd\sin\phi\cos\phi]$.\\
$<\DF{\partial h_1}{\partial x},\DF{\partial f}{\partial
y}>=N^2r^{2(N-1)}[ab\cos^2\phi-a^2\sin\phi\cos\phi+
b^2\sin\phi\cos\phi-ab\sin^2\phi+cd\cos^2\phi-c^2\sin\phi\cos\phi+d^2\sin\phi\cos\phi-cd\sin^2\phi].$\\
Using the conformal condition, and after simplification, we have
$$(a^2+c^2-b^2-d^2)\cos^2\phi+(b^2+d^2-a^2-c^2)\sin^2\phi+(4ab+4cd)\sin\phi\cos\phi=0,$$
and
$$(ab+cd)\cos^2\phi-(ab+cd)\sin^2\phi+(b^2+d^2-a^2-c^2)\sin\phi\cos\phi=0.$$
While the first one can be further reduced to
$$(a^2+c^2-b^2-d^2)\cos(2\phi)+(2ab+2cd)\sin(2\phi)=0,$$
and the second one can be reduced to
$$(ab+cd)\cos(2\phi)+\DF{b^2+d^2-a^2-c^2}{2}\sin(2\phi)=0.$$
In a matrix form, that is equivalent to
\[
\left(\begin{array}{cc}
a^2+c^2-b^2-d^2&2(ab+cd)\\
ab+cd&-\DF{a^2+c^2-b^2-d^2}{2}
\end{array}\right)
\left(\begin{array}{c}
\cos(2\phi)\\
\sin(2\phi)
\end{array}\right)=
\left(\begin{array}{c}
0\\
0
\end{array}\right)
\]
for any $\phi=(N-1)\theta$. Therefore, we must have
$$a^2+c^2-b^2-d^2=0, ab+cd=0.$$
Similarly, we have
$$\tilde{a}^2+\tilde{c}^2-\tilde{b}^2-\tilde{d}^2=0,\tilde{a}\tilde{b}+\tilde{c}\tilde{d}=0.$$
Now we will discuss the solutions of above equations.\\
If $c=0$, then $ab=0$, i.e $a=0$ or $b=0$. \\
If $c=a=0$, then $b=d=0$.\\
If $c=b=0$, then $a=\pm d$.\\
If $d=0$, then $ab=0$, i.e. $a=0$ or $b=0$.\\
If $d=a=0$, then $c=\pm b$.\\
If $d=b=0$, then $a=c=0$.\\
Now we assume that $cd\ne 0$, let $a=kc$, $b=ld$, for some constants $k,l$.\\
$$a^2+c^2-b^2-d^2=k^2c^2+c^2-l^2d^2-d^2=(k^2+1)c^2-(l^2+1)d^2=0$$
$$ab+cd=(kl+1)cd=0.$$
Since $cd\ne 0$, $kl=-1$, i.e. $l=-1/k$.\\
So $(k^2+1)c^2=(l^2+1)d^2=(\DF{1}{k^2}+1)d^2=\DF{k^2+1}{k^2}d^2$.\\
Hence $d^2=k^2c^2$, i.e. $d=\pm kc$, $b=ld=-\DF{1}{k}\cdot\pm kc=\mp c$.\\
In a word, here are the possible solutions of $a,b,c,d$:
\begin{equation*}
\begin{split}
(1)\;&(a,b,c,d)=(d,0,0,d),d\ne 0\\
(2)\;&(a,b,c,d)=(-d,0,0,d), d\ne 0\\
(3)\;&(a,b,c,d)=(0,b,b,0), b\ne 0\\
(4)\;&(a,b,c,d)=(0,b,-b,0), b\ne 0\\
(5)\;&(a,b,c,d)=(kc,-c,c,kc), k\ne 0, c\ne 0\\
(6)\;&(a,b,c,d)=(kc,c,c,-kc), k\ne 0, c\ne 0\\
(7)\;&(a,b,c,d)=(0,0,0,0)
\end{split}
\end{equation*}
We have the same conclusions about $\tilde{a},\tilde{b},\tilde{c},\tilde{d}$.\\
Finally, we will check the requirement that $\eta\circ
g=\DF{1}{2}(r^Ng_1(\theta)+r^Ng_2(\theta))$ is Dir minimizing, i.e
$r^N((a+\tilde{a})\cos(N\theta)+(b+\tilde{b})\sin(N\theta),(c+\tilde{c})\cos(N\theta)+(d+\tilde{d})\sin(N\theta))$
is Dir minimizing. Therefore,
the 4-tuple $(a+\tilde{a},b+\tilde{b},c+\tilde{c},d+\tilde{d})$ must be in one of the seven forms above.\\
Now let us consider the possibility of matching up the two 4-tuples
$(a,b,c,d),$
and $(\tilde{a},\tilde{b},\tilde{c},\tilde{d})$.\\
(1)+(1),i.e $(a,b,c,d)=(d,0,0,d), d\ne 0$,
$(\tilde{a},\tilde{b},\tilde{c},\tilde{d})=(\tilde{d},0,0,\tilde{d}),
\tilde{d}\ne 0$.
$$(a+\tilde{a},b+\tilde{b},c+\tilde{c},d+\tilde{d})=(d+\tilde{d},0,0,d+\tilde{d})$$
$$g(r,0)=[[r^N(d,0)]]+[[r^N(\tilde{d},0)]]$$
$$g(r,2\pi)=[[r^N(d\cos(2\pi N),d\sin(2\pi N))]]+[[r^N(\tilde{d}\cos(2\pi N),\tilde{d}\sin(2\pi N))]]$$
Let $\psi=2\pi N$. Since $g$ is H$\ddot{o}$lder continuous in
$\mathbb{U}_1^m(0)$, $g(r,0)=g(r,2\pi)$, i.e.,
$$d\sin(2\pi N)=\tilde{d}\sin(2\pi N)=0$$
$$\sin(2\pi N)=0$$
$$2\pi N=k\pi,\mbox{i.e}\;N=k/2, \mbox{for some integer}\;k=1,2,\cdot\cdot\cdot$$
Case 1:
$$d\cos(2\pi N)=d,\tilde{d}\cos(2\pi N)=\tilde{d}$$
$$\cos(2\pi N)=1,\mbox{i.e}\;2\pi N=2k\pi,\mbox{for some integer}\;k=1,2,\cdot\cdot\cdot$$
$$N=k,\mbox{for some integer}\;k=1,2,\cdot\cdot\cdot$$
Therefore,
$g(r,\theta)=dr^k[[(\cos(k\theta),\sin(k\theta))]]+\tilde{d}r^k[[(\cos(k\theta),\sin(k\theta))]]$.
\\
Case 2:
$$d=\tilde{d}\cos(2\pi N),\tilde{d}=d\cos(2\pi N)$$
$$\cos(2\pi N)=\pm 1,\mbox{i.e.} 2\pi N=k\pi$$
$$N=k/2,\mbox{for some integer}\;k=1,2,\cdot\cdot\cdot$$
If $N$ is an integer, then $d=\tilde{d}$.
$$g(r,\theta)=2dr^k[[(\cos(k\theta),\sin(k\theta))]],$$
in which case $\Sigma(g)=\emptyset$, a contradiction to our
assumption. Hence $N=k/2, \mbox{for some odd integer}\;k$
$$d=-\tilde{d}$$
$$g(r,\theta)=dr^{k/2}[[(\cos(\theta k/2),\sin(\theta k/2))]]+(-d)r^{k/2}[[(\cos(\theta k/2),\sin(\theta k/2))]]$$
(1)+(2): $(a,b,c,d)=(d,0,0,d),d\ne
0,(\tilde{a},\tilde{b},\tilde{c},\tilde{d})=(-\tilde{d},0,0,\tilde{d}),\tilde{d}\ne
0$.
$$(a,b,c,d)+(\tilde{a},\tilde{b},\tilde{c},\tilde{d})=(d-\tilde{d},0,0,d+\tilde{d})$$
which is in neither of those seven forms. \\
(1)+(3): $(a,b,c,d)=(d,0,0,d),d\ne 0,
(\tilde{a},\tilde{b},\tilde{c},\tilde{d})=(0,\tilde{b},\tilde{b},0),\tilde{d}\ne
0$
$$(a,b,c,d)+(\tilde{a},\tilde{b},\tilde{c},\tilde{d})=(d,\tilde{b},\tilde{b},d)$$
which is in neither of those seven forms.\\
(1)+(4): $(a,b,c,d)=(d,0,0,d),d\ne
0,(\tilde{a},\tilde{b},\tilde{c},\tilde{d})=(0,\tilde{b},-\tilde{b},0),\tilde{d}\ne
0$.
$$(a,b,c,d)+(\tilde{a},\tilde{b},\tilde{c},\tilde{d})=(d,\tilde{b},-\tilde{b},d)$$
which is of form (5).
$$g(r,0)=[[r^N(d,0)]]+[[r^N(0,-\tilde{b})]]$$
$$g(r,2\pi)=[[r^N(d\cos(2\pi N),d\sin(2\pi N))]]+[[r^N(\tilde{b}\sin(2\pi N),-\tilde{b}\cos(2\pi N))]]$$
Case 1:
$$d=d\cos(2\pi N),0=d\sin(2\pi N),0=\tilde{b}\sin(2\pi N),-\tilde{b}=-\tilde{b}\cos(2\pi N)$$
$$\cos(2\pi N)=1,\sin(2\pi N)=0$$
$$N=k,\mbox{for some integer}\;k=1,2,\cdot\cdot\cdot$$
$$g(r,\theta)=dr^k[[(\cos(k\theta),\sin(k\theta))]]+\tilde{b}r^k[[(\sin(k\theta),-\cos(k\theta))]]$$
 Case 2:
$$d=\tilde{b}\sin(2\pi N),0=-\tilde{b}\cos(2\pi N), 0=d\cos(2\pi N), -\tilde{b}=d\sin(2\pi N)$$
which has no solutions. \\
(1)+(5). $(a,b,c,d)=(d,0,0,d),d\ne
0,(\tilde{a},\tilde{b},\tilde{c},\tilde{d})=(l\tilde{c},-\tilde{c},\tilde{c},l\tilde{c}),l\ne
0, \tilde{c}\ne 0$.
$$(a,b,c,d)+(\tilde{a},\tilde{b},\tilde{c},\tilde{d})=(d+l\tilde{c},-\tilde{c},\tilde{c},d+l\tilde{c})$$
$$g(r,0)=[[r^N(d,0)]]+[[r^N(l\tilde{c},\tilde{c})]]$$
$$g(r,2\pi)=[[r^N(d\cos(2\pi N),d\sin(2\pi N))]]+$$
$$[[r^N(l\tilde{c}\cos(2\pi N)-\tilde{c}\sin(2\pi N),\tilde{c}
\cos(2\pi N)+l\tilde{c}\sin(2\pi N))]]$$ Case 1:
$$d=l\tilde{c}\cos\psi-\tilde{c}\sin\psi,0=\tilde{c}\cos\psi+l\tilde{c}\sin\psi$$
$$l\tilde{c}=d\cos\psi,\tilde{c}=d\sin\psi$$
No solution.\\
Case 2:
$$d=d\cos\psi,0=d\sin\psi,$$
$$l\tilde{c}=l\tilde{c}\cos\psi-\tilde{c}\sin\psi,\tilde{c}=\tilde{c}\cos\psi+l\tilde{c}\sin\psi$$
The solution is $\cos\psi=1,\sin\psi=0$, hence $2\pi N=\psi=2k\pi$,
$N=k$.
$$g(r,\theta)=dr^k[[(\cos(k\theta),\sin(k\theta))]]+\tilde{c}r^k[[(l\cos(k\theta)-\sin(k\theta),\cos(k\theta)+l\sin(k\theta))]]$$
(1)+(6). $(a,b,c,d)=(d,0,0,d),d\ne 0,
(\tilde{a},\tilde{b},\tilde{c},\tilde{d})=(l\tilde{c},\tilde{c},\tilde{c},-l\tilde{c}),
l\ne 0,\tilde{c}\ne 0$.\\
$$(a+\tilde{a},b+\tilde{b},c+\tilde{c},d+\tilde{d})=(d+l\tilde{c},\tilde{c},\tilde{c},d-l\tilde{c})$$
which is in neither of the seven forms above.\\
(1)+(7). $(a,b,c,d)=(d,0,0,d),d\ne 0,
(\tilde{a},\tilde{b},\tilde{c},\tilde{d})=(0,0,0,0)$
$$(a,b,c,d)+(\tilde{a},\tilde{b},\tilde{c},\tilde{d})=(d,0,0,d)$$
$$g(r,0)=r^N[[(d,0)]]+r^N[[(0,0)]]$$
$$g(r,2\pi)=r^N[[(d\cos(2\pi N),d\sin(2\pi N))]]+r^N[[(0,0)]]$$
Since $d\ne 0$, the only possible matching up is
$$d=d\cos\psi,0=d\sin\psi$$
i.e $\sin\psi=0,\cos\psi=1$.
$$N=k,\;\mbox{for some positivie integer}\;k=1,2,\cdot\cdot\cdot$$
$$g(r,\theta)=dr^k[[(\cos(k\theta),\sin(k\theta))]]+[[(0,0)]]$$
(2)+(2). $(a,b,c,d)=(-d,0,0,d),d\ne
0,(\tilde{a},\tilde{b},\tilde{c},\tilde{d})=(-\tilde{d},0,0,\tilde{d}),\tilde{d}\ne
0$.
$$(a,b,c,d)+(\tilde{a},\tilde{b},\tilde{c},\tilde{d})=(-d-\tilde{d},0,0,d+\tilde{d})$$
$$g(r,0)=[[r^N(-d,0)]]+[[r^N(-\tilde{d},0)]]$$
$$g(r,2\pi)=[[r^N(-d\cos\psi,d\sin\psi)]]+[[r^N(-\tilde{d}\cos\psi,\tilde{d}\sin\psi)]]$$
Case 1:
$$-d=-d\cos\psi,0=d\sin\psi,-\tilde{d}=-\tilde{d}\cos\psi,0=\tilde{d}\sin\psi$$
Hence $\cos\psi=1,\sin\psi=0$, i.e.,
$$N=k,\;\mbox{for some positive integer}\;k=1,2,\cdot\cdot\cdot.$$
$$g(r,\theta)=dr^k[[(-\cos(k\theta),\sin(k\theta))]]+\tilde{d}r^k[[(-\cos(k\theta),\sin(k\theta))]]$$
Case 2:
$$-d=-\tilde{d}\cos\psi,0=\tilde{d}\sin\psi,-\tilde{d}=-d\cos\psi,0=d\sin\psi$$
The solution is $\sin\psi=0,\cos\psi=\pm 1$.\\
If $\cos\psi=1$, then $d=\tilde{d}$, which means that
$\Sigma(g)=\emptyset$, a contradiction to our assumption. Therefore
$\cos\psi=-1$, which means $N=k/2$ for some odd integer $k$.
Moreover, we get $d=-\tilde{d}$.
$$g(r,\theta)=\pm dr^k[[(-\cos(k\theta),\sin(k\theta))]]$$
(2)+(3). $(a,b,c,d)=(-d,0,0,d),d\ne
0,(\tilde{a},\tilde{b},\tilde{c},\tilde{d})=(0,\tilde{b},\tilde{b},0),\tilde{b}\ne
0$.
$$(a,b,c,d)+(\tilde{a},\tilde{b},\tilde{c},\tilde{d})=(-d,\tilde{b},\tilde{b},d)$$
$$g(r,0)=[[r^N(-d,0)]]+[[r^N(0,\tilde{b})]]$$
$$g(r,2\pi)=[[r^N(-d\cos\psi,d\sin\psi)]]+[[r^N(\tilde{b}\sin\psi,\tilde{b}\cos\psi)]]$$
Case 1:
$$-d=-d\cos\psi,0=d\sin\psi,0=\tilde{b}\sin\psi,\tilde{b}=\tilde{b}\cos\psi$$
Hence $\sin\psi=0,\cos\psi=1$, i.e. $N=k,\;\mbox{for some positive
integer}\; k=1,2,\cdot\cdot\cdot$.
$$g(r,\theta)=d[[r^k(-\cos(k\theta),\sin(k\theta))]]+\tilde{b}r^k[[(\sin(k\theta),\cos(k\theta))]]$$
Case 2:
$$-d=\tilde{b}\sin\psi,0=\tilde{b}\cos\psi,0=-d\cos\psi,\tilde{b}=d\sin\psi$$
No solution.\\
(2)+(4) $(a,b,c,d)=(-d,0,0,d),d\ne
0,(\tilde{a},\tilde{b},\tilde{c},\tilde{d})=(0,\tilde{b},-\tilde{b},0),\tilde{b}\ne
0$.
$$(a,b,c,d)+(\tilde{a},\tilde{b},\tilde{c},\tilde{d})=(-d,\tilde{b},-\tilde{b},d)$$
which is in neither of the seven forms above.\\
(2)+(5) $(a,b,c,d)=(-d,0,0,d),d\ne
0,(\tilde{a},\tilde{b},\tilde{c},\tilde{d})=(l\tilde{c},-\tilde{c},\tilde{c},l\tilde{c}),
l\ne 0,\tilde{c}\ne 0$
$$(a,b,c,d)+(\tilde{a},\tilde{b},\tilde{c},\tilde{d})=(l\tilde{c}-d,-\tilde{c},\tilde{c},l\tilde{c}+d)$$
which is in neither of the seven forms above.\\
(2)+(6) $(a,b,c,d)=(-d,0,0,d),d\ne
0,(\tilde{a},\tilde{b},\tilde{c},\tilde{d})=(l\tilde{c},\tilde{c},\tilde{c},-l\tilde{c}),
l\ne 0,\tilde{c}\ne 0$
$$(a,b,c,d)+(\tilde{a},\tilde{b},\tilde{c},\tilde{d})=(l\tilde{c}-d,\tilde{c},\tilde{c},d-l\tilde{c})$$
$$g(r,0)=r^N[[(-d,0)]]+r^N[[(l\tilde{c},\tilde{c})]]$$
$$g(r,2\pi)=r^N[[(-d\cos\psi,d\sin\psi)]]+r^N[[(l\tilde{c}\cos\psi+\tilde{c}\sin\psi,\tilde{c}\cos\psi-l\tilde{c}\sin\psi)]]$$
Case 1:
$$-d=-d\cos\psi,0=d\sin\psi$$
$$l\tilde{c}=l\tilde{c}\cos\psi+\tilde{c}\sin\psi,\tilde{c}=\tilde{c}\cos\psi-l\tilde{c}\sin\psi$$
The solution is $\sin\psi=0,\cos\psi=1$, i.e. $N=k,\;\mbox{for some
potitive integer}\;k$.
$$g(r,\theta)=dr^k[[(-\cos(k\theta),\sin(k\theta))]]+\tilde{c}r^k[[(l\cos(k\theta)+\sin(k\theta),\cos(k\theta)-l\sin(k\theta))]]$$
 Case 2:
$$-d=l\tilde{c}\cos\psi+\tilde{c}\sin\psi,0=\tilde{c}\cos\psi-l\tilde{c}\sin\psi$$
$$l\tilde{c}=-d\cos\psi,\tilde{c}=d\sin\psi$$
No solution.\\
(2)+(7). $(a,b,c,d)=(-d,0,0,d),d\ne
0,(\tilde{a},\tilde{b},\tilde{c},\tilde{d})=(0,0,0,0)$
$$(a,b,c,d)+(\tilde{a},\tilde{b},\tilde{c},\tilde{d})=(-d,0,0,d)$$
$$g(r,0)=r^N[[(-d,0)]]+r^N[[(0,0)]]$$
$$g(r,2\pi)=r^N[[(-d\cos\psi,d\sin\psi)]]+r^N[[(0,0)]]$$
Since $d\ne 0$, there is only one way of matching up:
$$-d=-d\cos\psi,0=d\sin\psi$$
Therefore $N=k$ for some positive integer $k$.
$$g(r,\theta)=r^k[[(-d\cos(k\theta),d\sin(k\theta))]]+[[(0,0)]]$$
(3)+(3): $(a,b,c,d)=(0,b,b,0),b\ne
0,(\tilde{a},\tilde{b},\tilde{c},\tilde{d})=(0,\tilde{b},\tilde{b},0),\tilde{b}\ne
0$
$$(a,b,c,d)+(\tilde{a},\tilde{b},\tilde{c},\tilde{d})=(0,b+\tilde{b},b+\tilde{b},0)$$
$$g(r,0)=r^N[[(0,b)]]+r^N[[(0,\tilde{b})]]$$
$$g(r,2\pi)=r^N[[(b\sin\psi,b\cos\psi)]]+r^N[[(\tilde{b}\sin\psi,\tilde{b}\cos\psi)]]$$
Case 1:
$$0=b\sin\psi,b=b\cos\psi,0=\tilde{b}\sin\psi,\tilde{b}=\tilde{b}\cos\psi$$
The solution is $\sin\psi=0,\cos\psi=1$, i.e. $N=k\;\mbox{for some
positive integer}\;k$.
$$g(r,\theta)=br^k[[(\sin(k\theta),\cos(k\theta))]]+\tilde{b}r^k[[(\sin(k\theta),\cos(k\theta))]]$$
Case 2:
$$0=\tilde{b}\sin\psi,b=\tilde{b}\cos\psi,0=b\sin\psi,\tilde{b}=b\cos\psi$$
The solution is $\sin\psi=0,\cos\psi=\pm 1$. \\
If $\cos\psi=1$, then $b=\tilde{b}$, which means
$\Sigma(g)=\emptyset$, a contradiction to our assumption.Therefore
$\cos\psi=-1$ ,$b=-\tilde{b}$ and $N=k/2$ for some odd integer $k$.
$$g(r,\theta)=br^{k/2}[[(\sin(\theta k/2),\cos(\theta k/2))]]+(-b)r^{k/2}[[(\sin(\theta k/2),\cos(\theta k/2))]]$$
(3)+(4) $(a,b,c,d)=(0,b,b,0),b\ne 0,
(\tilde{a},\tilde{b},\tilde{c},\tilde{d})=(0,\tilde{b},-\tilde{b},0),
\tilde{b}\ne 0$.
$$(a,b,c,d)+(\tilde{a},\tilde{b},\tilde{c},\tilde{d})=(0,b+\tilde{b},b-\tilde{b},0)$$
which is in neither of the seven forms above.\\
(3)+(5) $(a,b,c,d)=(0,b,b,0),b\ne 0,
(\tilde{a},\tilde{b},\tilde{c},\tilde{d})=(l\tilde{c},-\tilde{c},\tilde{c},l\tilde{c}),
l\ne 0,\tilde{c}\ne 0$.
$$(a,b,c,d)+(\tilde{a},\tilde{b},\tilde{c},\tilde{d})=(l\tilde{c},b-\tilde{c},b+\tilde{c},l\tilde{c})$$
which is in neither of the seven forms above.\\
(3)+(6) $(a,b,c,d)=(0,b,b,0),b\ne 0,
(\tilde{a},\tilde{b},\tilde{c},\tilde{d})=(l\tilde{c},\tilde{c},\tilde{c},-l\tilde{c}),
l\ne 0,\tilde{c}\ne 0$.
$$(a,b,c,d)+(\tilde{a},\tilde{b},\tilde{c},\tilde{d})=(l\tilde{c},b+\tilde{c},b+\tilde{c},-l\tilde{c})$$
$$g(r,0)=r^N[[(0,b)]]+r^N[[(l\tilde{c},\tilde{c})]]$$
$$g(r,2\pi)=r^N[[(b\sin\psi,b\cos\psi)]]+r^N[[(l\tilde{c}\cos\psi+\tilde{c}\sin\psi,\tilde{c}\cos\psi-l\tilde{c}\sin\psi)]]$$
Case 1:
$$0=b\sin\psi,b=b\cos\psi$$
$$l\tilde{c}=l\tilde{c}\cos\psi+\tilde{c}\sin\psi,\tilde{c}=\tilde{c}\cos\psi-l\tilde{c}\sin\psi$$
The solution is $\sin\psi=0,\cos\psi=1$, hence $N=k\;\mbox{for some
positive integer}\;k$.
$$g(r,\theta)=br^k[[(\sin(k\theta),\cos(k\theta))]]+r^k[[(l\tilde{c}\cos(k\theta)+\tilde{c}\sin(k\theta),\tilde{c}\cos(k\theta)-l\tilde{c}\sin(k\theta))]]$$
Case 2:
$$0=l\tilde{c}\cos\psi+\tilde{c}\sin\psi,b=\tilde{c}\cos\psi-l\tilde{c}\sin\psi$$
$$l\tilde{c}=b\sin\psi,\tilde{c}=b\cos\psi$$
No solution.\\
(3)+(7) $(a,b,c,d)=(0,b,b,0),b\ne 0,
(\tilde{a},\tilde{b},\tilde{c},\tilde{d})=(0,0,0,0)$.
$$(a,b,c,d)+(\tilde{a},\tilde{b},\tilde{c},\tilde{d})=(0,b,b,0)$$
$$g(r,0)=r^N[[(0,b)]]+[[(0,0)]]$$
$$g(r,2\pi)=r^N[[(b\sin\psi,b\cos\psi)]]+[[(0,0)]]$$
Since $b\ne 0$, there is only one way of matching up:
$$0=b\sin\psi,b=b\cos\psi$$
Hence $\sin\psi=0,\cos\psi=1$, i.e. $N=k$ for some positive integer
$k$.
$$g(r,\theta)=br^k[[(\sin(k\theta),\cos(k\theta))]]+[[(0,0)]]$$
(4)+(4) $(a,b,c,d)=(0,b,-b,0),b\ne 0,
(\tilde{a},\tilde{b},\tilde{c},\tilde{d})=(0,\tilde{b},-\tilde{b},0),\tilde{b}\ne
0$.
$$(a,b,c,d)+(\tilde{a},\tilde{b},\tilde{c},\tilde{d})=(0,b+\tilde{b},-(b+\tilde{b}),0)$$
$$g(r,0)=r^N[[(0,-b)]]+r^N[[(0,-\tilde{b})]]$$
$$g(r,2\pi)=r^N[[(b\sin\psi,-b\cos\psi)]]+r^N[[(\tilde{b}\sin\psi,-\tilde{b}\cos\psi)]]$$
Case 1:
$$0=b\sin\psi,-b=-b\cos\psi,0=\tilde{b}\sin\psi,-\tilde{b}=-\tilde{b}\cos\psi$$
The solution is $N=k\;\mbox{for some positive integer}\;k$.
$$g(r,\theta)=br^k[[(\sin(k\theta),-\cos(k\theta))]]+\tilde{b}r^k[[(\sin(k\theta),-\cos(k\theta))]]$$
Case 2:
$$0=\tilde{b}\sin\psi,-b=-\tilde{b}\cos\psi,0=b\sin\psi,-\tilde{b}=-b\cos\psi$$
The solution is $\sin\psi=0,\cos\psi=\pm 1$.\\
If $\cos\psi=1$, then $b=\tilde{b}$, which means
$\Sigma(g)=\emptyset$, a contradiction to our assumption. Therefore
$\cos\psi=-1$, hence $b=-\tilde{b}$, and $N=k/2$ for some odd
integer $k$.
$$g(r,\theta)=br^{k/2}[[(\sin(\theta k/2),-\cos(\theta k/2))]]+(-b)r^{k/2}[[(\sin(\theta k/2),-\cos(\theta k/2))]]$$
(4)+(5) $(a,b,c,d)=(0,b,-b,0),b\ne 0,
(\tilde{a},\tilde{b},\tilde{c},\tilde{d})=(l\tilde{c},-\tilde{c},\tilde{c},l\tilde{c}),
l\ne 0,\tilde{c}\ne 0$
$$(a,b,c,d)+(\tilde{a},\tilde{b},\tilde{c},\tilde{d})=(l\tilde{c},b-\tilde{c},-b+\tilde{c},l\tilde{c})$$
$$g(r,0)=r^N[[(0,-b)]]+r^N[[(l\tilde{c},\tilde{c})]]$$
$$g(r,2\pi)=r^N[[(b\sin\psi,-b\cos\psi)]]+r^N[[(l\tilde{c}\cos\psi-\tilde{c}\sin\psi,\tilde{c}\cos\psi+l\tilde{c}\sin\psi)]]$$
Case 1
$$0=b\sin\psi,-b=-b\cos\psi$$
$$l\tilde{c}=l\tilde{c}\cos\psi-\tilde{c}\sin\psi,\tilde{c}=\tilde{c}\cos\psi+l\tilde{c}\sin\psi$$
The solution is $N=k\;\mbox{for some positive integer}\;k$.
$$g(r,\theta)=br^k[[(\sin(k\theta),-\cos(k\theta))]]+\tilde{c}r^k[[(l\cos(k\theta)-\sin(k\theta),\cos(k\theta)+l\sin(k\theta))]]$$
Case 2
$$0=l\tilde{c}\cos\psi-\tilde{c}\sin\psi,-b=\tilde{c}\cos\psi+l\tilde{c}\sin\psi$$
$$l\tilde{c}=b\sin\psi,\tilde{c}=-b\cos\psi$$
No solutions.\\
(4)+(6) $(a,b,c,d)=(0,b,-b,0),b\ne
0,(\tilde{a},\tilde{b},\tilde{c},\tilde{d})=(l\tilde{c},\tilde{c},\tilde{c},-l\tilde{c}),
l\ne 0,\tilde{c}\ne 0$.
$$(a,b,c,d)+(\tilde{a},\tilde{b},\tilde{c},\tilde{d})=(l\tilde{c},b+\tilde{c},-b+\tilde{c},-l\tilde{c})$$
which is in neither of the seven forms above.\\
(4)+(7) $(a,b,c,d)=(0,b,-b,0),b\ne
0,(\tilde{a},\tilde{b},\tilde{c},\tilde{d})=(0,0,0,0)$.
$$(a,b,c,d)+(\tilde{a},\tilde{b},\tilde{c},\tilde{d})=(0,b,-b,0)$$
$$g(r,0)=r^N[[(0,-b)]]+[[(0,0)]]$$
$$g(r,2\pi)=r^N[[(b\sin\psi,-b\cos\psi)]]+[[(0,0)]]$$
Since $b\ne 0$, there is only one way of matching up
$$0=b\sin\psi,-b=-b\cos\psi$$
Hence $N=k$ for some positive integer $k$.
$$g(r,\theta)=r^k[[(b\sin(k\theta),-b\cos(k\theta))]]+[[(0,0)]]$$
(5)+(5) $(a,b,c,d)=(lc,-c,c,lc),l\ne 0, c\ne
0,(\tilde{a},\tilde{b},\tilde{c},\tilde{d})=(\tilde{l}\tilde{c},-\tilde{c},\tilde{c},\tilde{l}\tilde{c}),
\tilde{l}\ne 0,\tilde{c}\ne 0$
$$(a+\tilde{a},b+\tilde{b},c+\tilde{c},d+\tilde{d})=(lc+\tilde{l}\tilde{c},-c-\tilde{c},c+\tilde{c},lc+\tilde{l}\tilde{c})$$
$$g(r,0)=r^N[[(lc,c)]]+r^N[[(\tilde{l}\tilde{c},\tilde{c})]]$$
$$g(r,2\pi)=r^N[[(lc\cos\psi-c\sin\psi,c\cos\psi+lc\sin\psi)]]+r^N[[(\tilde{l}\tilde{c}\cos\psi-\tilde{c}\sin\psi,\tilde{c}\cos\psi+
\tilde{l}\tilde{c}\sin\psi)]]$$ Case 1
$$lc=lc\cos\psi-c\sin\psi,c=c\cos\psi+lc\sin\psi$$
$$\tilde{l}\tilde{c}=\tilde{l}\tilde{c}\cos\psi-\tilde{c}\sin\psi,\tilde{c}=\tilde{c}\cos\psi+\tilde{l}\tilde{c}\sin\psi$$
The solution is $N=k\;\mbox{for some positive integer}\;k$.
$$g(r,\theta)=r^k[[(lc\cos(k\theta)-c\sin(k\theta),c\cos(k\theta)+lc\sin(k\theta))]]$$
$$+
r^k[[(\tilde{l}\tilde{c}\cos(k\theta)-\tilde{c}\sin(k\theta),\tilde{c}\cos(k\theta)+\tilde{l}\tilde{c}\sin(k\theta))]]$$
Case 2
$$lc=\tilde{l}\tilde{c}\cos\psi-\tilde{c}\sin\psi,c=\tilde{c}\cos\psi+\tilde{l}\tilde{c}\sin\psi$$
$$\tilde{l}\tilde{c}=lc\cos\psi-c\sin\psi,\tilde{c}=c\cos\psi+lc\sin\psi$$
The solution is $\sin\psi=0,\cos\psi=\pm 1$. \\
If $\cos\psi=1$, then $l=\tilde{l}$ and $c=\tilde{c}$, which means
$\Sigma(g)=\emptyset$, a contradiction to our assumption. Therefore
$\cos\psi=-1$, hence $N=k/2,\;\mbox{for some odd integer}\;k$,
$c=-\tilde{c},l=\tilde{l}$
$$g(r,\theta)=\pm r^{k/2}[[(lc\cos(\theta k/2)-c\sin(\theta k/2),c\cos(\theta k/2)+lc\sin(\theta k/2))]]$$
(5)+(6) $(a,b,c,d)=(lc,-c,c,lc),l\ne 0, c\ne 0,
(\tilde{a},\tilde{b},\tilde{c},\tilde{d})=(\tilde{l}\tilde{c},\tilde{c},\tilde{c},-\tilde{l}\tilde{c}),
\tilde{l}\ne 0,\tilde{c}\ne 0$
$$(a,b,c,d)+(\tilde{a},\tilde{b},\tilde{c},\tilde{d})=(lc+\tilde{l}\tilde{c},-c+\tilde{c},c+\tilde{c},lc-\tilde{l}\tilde{c})$$
which is in neither of the seven forms above.\\
(5)+(7) $(a,b,c,d)=(lc,-c,c,lc),l\ne 0,c\ne
0,(\tilde{a},\tilde{b},\tilde{c},\tilde{d})=(0,0,0,0)$.
$$(a,b,c,d)+(\tilde{a},\tilde{b},\tilde{c},\tilde{d})=(lc,-c,c,lc)$$
$$g(r,0)=r^N[[(lc,c)]]+[[(0,0)]]$$
$$g(r,2\pi)=r^N[[(lc\cos\psi-c\sin\psi,c\cos\psi+lc\sin\psi)]]+[[(0,0)]]$$
Since $c\ne 0$, there is only one way of matching up
$$lc=lc\cos\psi-c\sin\psi,c=c\cos\psi+lc\sin\psi$$
The solution is $\sin\psi=0,\cos\psi=1$, hence $N=k$ for some
positive integer $k$.
$$g(r,\theta)=r^k[[(lc\cos(k\theta)-c\sin(k\theta),c\cos(k\theta)+lc\sin(k\theta))]]+[[(0,0)]]$$
(6)+(6) $(a,b,c,d)=(lc,c,c,-lc), l\ne 0,c\ne 0,
(\tilde{a},\tilde{b},\tilde{c},\tilde{c})=(\tilde{l}\tilde{c},\tilde{c},\tilde{c},-\tilde{l}\tilde{c}),
\tilde{l}\ne 0,\tilde{c}\ne 0$
$$(a,b,c,d)+(\tilde{a},\tilde{b},\tilde{c},\tilde{d})=(lc+\tilde{l}\tilde{c},c+\tilde{c},c+\tilde{c},-lc-\tilde{l}\tilde{c})$$
$$g(r,0)=r^N[[(lc,c)]]+r^N[[(\tilde{l}\tilde{c},\tilde{c})]]$$
$$g(r,2\pi)=r^N[[(lc\cos\psi+c\sin\psi,c\cos\psi-lc\sin\psi)]]+r^N[[(\tilde{l}\tilde{c}\cos\psi+\tilde{c}\sin\psi,
\tilde{c}\cos\psi-\tilde{l}\tilde{c}\sin\psi)]]$$ Case 1:
$$lc=lc\cos\psi+c\sin\psi,c=c\cos\psi-lc\sin\psi$$
$$\tilde{l}\tilde{c}=\tilde{l}\tilde{c}\cos\psi+\tilde{c}\sin\psi,\tilde{c}=\tilde{c}\cos\psi-\tilde{l}\tilde{c}\sin\psi$$
The solution is $\sin\psi=0,\cos\psi=1$, i.e $N=k\;\mbox{for some
positive integer}\;k$.
$$g(r,\theta)=r^k[[(lc\cos(k\theta)+c\sin(k\theta),c\cos(k\theta)-lc\sin(k\theta))]]$$
$$+r^k[[(\tilde{l}\tilde{c}\cos(k\theta)
+\tilde{c}\sin(k\theta),\tilde{c}\cos(k\theta)-\tilde{l}\tilde{c}\sin(k\theta))]]$$
 Case 2:
$$lc=\tilde{l}\tilde{c}\cos\psi+\tilde{c}\sin\psi,c=\tilde{c}\cos\psi-\tilde{l}\tilde{c}\sin\psi$$
$$\tilde{l}\tilde{c}=lc\cos\psi+c\sin\psi,\tilde{c}=c\cos\psi-lc\sin\psi$$
The solution is $\sin\psi=0,\cos\psi=\pm 1$. \\
If $\cos\psi=1$, then $c=\tilde{c},l=\tilde{l}$, which means
$\Sigma(g)=\emptyset$, a contradiction to our assumption. Therefore
$\cos\psi=-1$, hence $N=k/2\;\mbox{for some odd integer}\;k$ and
$c=-\tilde{c},l=\tilde{l}$.
$$g(r,\theta)=\pm r^{k/2}[[(lc\cos(\theta k/2)+c\sin(\theta k/2),c\cos(\theta k/2)-lc\sin(\theta k/2))]]$$
(6)+(7) $(a,b,c,d)=(lc,c,c,-lc), l\ne 0,c\ne
0,(\tilde{a},\tilde{b},\tilde{c},\tilde{d})=(0,0,0,0)$
$$(a,b,c,d)+(\tilde{a},\tilde{b},\tilde{c},\tilde{c})=(lc,c,c,-lc)$$
$$g(r,0)=r^N[[(lc,c)]]+[[(0,0)]]$$
$$g(r,2\pi)=r^N[[(lc\cos\psi+c\sin\psi,c\cos\psi-lc\sin\psi)]]+[[(0,0)]]$$
Since $c\ne 0$, the only matching up is
$$lc=lc\cos\psi+c\sin\psi,c=c\cos\psi-lc\sin\psi$$
The solution is $\sin\psi=0,\cos\psi=1$, hence $N=k\;\mbox{for some
positive integer}\;k$.
$$g(r,\theta)=r^k[[(lc\cos(k\theta)+c\sin(k\theta),c\cos(k\theta)-lc\sin(k\theta))]]+[[(0,0)]]$$
(7)+(7) This case does not happen because otherwise $g\equiv
2[[0]]$.
\subsection{$\Sigma=\emptyset$}
If $\Sigma=\emptyset$, applying Theorem 2.3(2) to the function
$g$, more specifically, by choosing $W=\mathbb{U}_1^2(0)$, we get
$$g(x)=2[[g_1(x)]],x\in \mathbb{U}_1^2(0),$$
for some minimizing harmonic function
$g_1:\mathbb{U}_1^2(0)\rightarrow
\mathbb{R}^2$.\\
The same argument as above tells that
$$g_1(r,\theta)=(r^N(a\cos(N\theta)+b\sin(N\theta)),r^N(c\cos(N\theta)+d\sin(N\theta)))$$
for some constants $a,b,c,d$ in one of the six forms (1)-(6) above.\\
Case (1). $(a,b,c,d)=(d,0,0,d)$.\\
$$g_1(r,0)=(ar^N,cr^N)=(dr^N,0)$$
\begin{equation*}
\begin{split}
g_1(r,2\pi)&=(r^N(a\cos\psi+b\sin\psi),r^N(c\cos\psi+d\sin\psi))\\
&=(r^N(d\cos\psi),r^N(d\sin\psi))
\end{split}
\end{equation*}
Therefore, $d=d\cos\psi,0=d\sin\psi$, i.e.,
$$2\pi N=2k\pi,\;\mbox{for some positive integer}\;k.
$$
Hence $N=k\;\mbox{for some positive integer}\; k$.\\
$$g(r,\theta)=2[[(dr^k\cos(k\theta),dr^k\sin(k\theta))]],$$
$\;\mbox{for some nonzero constant}\;d\;\mbox{and some positive
integer}
\;k.$\\
Case (2). $(a,b,c,d)=(-d,0,0,d)$.\\
$$g_1(r,0)=(ar^N,cr^N)=(-dr^N,0)$$
\begin{equation*}
\begin{split}
g_1(r,2\pi)&=(r^N(a\cos\psi+b\sin\psi),r^N(c\cos\psi+d\sin\psi))\\
&=(r^N(-d\cos\psi),r^N(d\sin\psi))
\end{split}
\end{equation*}
Therefore $-d=-d\cos\psi,0=d\sin\psi$, i.e.,
$$2\pi N=2k\pi,\;\mbox{for some positive integer}\;k$$
Hence $N=k\;\mbox{for some positive integer}\; k$.\\
$$g(r,\theta)=2[[(-dr^k\cos(k\theta),dr^k\sin(k\theta))]],$$
$\;\mbox{for some nonzero constant}\;d\;\mbox{and some positive
integer}
\;k.$\\
Case (3). $(a,b,c,d)=(0,b,b,0)$.\\
$$g_1(r,0)=(ar^N,cr^N)=(0,br^N)$$
\begin{equation*}
\begin{split}
g_1(r,2\pi)&=(r^N(a\cos\psi+b\sin\psi),r^N(c\cos\psi+d\sin\psi))\\
&=(r^N(b\sin\psi),r^N(b\cos\psi))
\end{split}
\end{equation*}
Therefore $0=b\sin\psi,b=b\cos\psi$, i.e.,
$$2\pi N=2k\pi,\;\mbox{for some positive integer}\;k$$
Hence $N=k\;\mbox{for some positive integer}\; k$.\\
$$g(r,\theta)=2[[(br^k\sin(k\theta),br^k\cos(k\theta))]],$$
$\;\mbox{for some nonzero constant}\;b\;\mbox{and some positive
integer}
\;k.$\\
Case (4). $(a,b,c,d)=(0,b,-b,0)$.\\
$$g_1(r,0)=(ar^N,cr^N)=(0,-br^N)$$
\begin{equation*}
\begin{split}
g_1(r,2\pi)&=(r^N(a\cos\psi+b\sin\psi),r^N(c\cos\psi+d\sin\psi))\\
&=(r^N(b\sin\psi),r^N(-b\cos\psi))
\end{split}
\end{equation*}
Therefore $0=b\sin\psi,-b=-b\cos\psi$, i.e.,
$$2\pi N=2k\pi,\;\mbox{for some positive integer}\;k$$
Hence $N=k\;\mbox{for some positive integer}\; k$.\\
$$g(r,\theta)=2[[(br^k\sin(k\theta),-br^k\cos(k\theta))]],$$
$\;\mbox{for some nonzero constant}\;b\;\mbox{and some positive
integer}
\;k.$\\
Case (5). $(a,b,c,d)=(lc,-c,c,lc)$.\\
$$g_1(r,0)=(ar^N,cr^N)=(lcr^N,cr^N)$$
\begin{equation*}
\begin{split}
g_1(r,2\pi)&=(r^N(a\cos\psi+b\sin\psi),r^N(c\cos\psi+d\sin\psi))\\
&=(r^N(lc\cos\psi-c\sin\psi),r^N(c\cos\psi+lc\sin\psi))
\end{split}
\end{equation*}
Therefore $lc=lc\cos\psi-c\sin\psi,c=c\cos\psi+lc\sin\psi$. Solving
that gives $\cos\psi=1,\sin\psi=0$, i.e.,
$$2\pi N=2k\pi,\;\mbox{for some positive integer}\;k$$
Hence $N=k\;\mbox{for some positive integer}\; k$.\\
$$g(r,\theta)=2[[(r^k(lc\cos(k\theta)-c\sin(k\theta)),r^k(c\cos(k\theta)+lc\sin(l\theta)))]],$$
$\;\mbox{for some nonzero constant}\;l,c\;\mbox{and some positive
integer}
\;k.$\\
Case (6). $(a,b,c,d)=(lc,c,c,-lc)$.\\
$$g_1(r,0)=(ar^N,cr^N)=(lcr^N,cr^N)$$
\begin{equation*}
\begin{split}
g_1(r,2\pi)&=(r^N(a\cos\psi+b\sin\psi),r^N(c\cos\psi+d\sin\psi))\\
&=(r^N(lc\cos\psi+c\sin\psi),r^N(c\cos\psi-lc\sin\psi))
\end{split}
\end{equation*}
Therefore $lc=lc\cos\psi+c\sin\psi,c=c\cos\psi-lc\sin\psi$. Solving
that gives $\cos\psi=1,\sin\psi=0$, i.e.,
$$2\pi N=2k\pi,\;\mbox{for some positive integer}\;k$$
Hence $N=k\;\mbox{for some positive integer}\; k$.\\
$$g(r,\theta)=2[[(r^k(lc\cos(k\theta)+c\sin(k\theta)),r^k(c\cos(k\theta)-lc\sin(k\theta)))]],$$
$\;\mbox{for some nonzero constant}\;l,c\;\mbox{and some positive
integer}
\;k.$\\
\end{proof}

\end{document}